\documentclass[12pt,a4paper]{article}
\usepackage{amsmath}
\usepackage{amsfonts}
\usepackage{amssymb}
\usepackage[dvipdfm]{graphicx}
\usepackage{graphics}
\usepackage{graphicx}
\usepackage[plainpages=false,pdfpagelabels,pagebackref=false,naturalnames=true,hyperindex=true,pdftitle={Sloane's Gap},pdfauthor={Hector Zenil}]{hyperref}

\begin{document}

\title{Sloane's Gap: Do Mathematical and Social Factors Explain the Distribution of Numbers in the OEIS?}
\author{Nicolas Gauvrit$^1$, Jean-Paul Delahaye$^2$ and Hector Zenil$^2$\\ 
$^{1}$LADR, EA 1547, Centre Chevaleret, Universit\'e de Paris VII\\adems@free.fr\\
$^{2}$LIFL, Laboratoire d'Informatique Fondamentale de Lille\\
UMR CNRS 8022, Universit\'e de Lille I\\
\{hector.zenil,delahaye\}@lifl.fr}

\maketitle

\begin{abstract}
The Online Encyclopedia of Integer Sequences (OEIS) is a catalog of integer sequences. We are particularly interested in the number of occurrences of $N(n)$ of an integer $n$ in the database. This number $N(n)$ marks the importance of $n$ and it varies noticeably from one number to another, and from one number to the next in a series. ``Importance" can be mathematically objective ($2^{10}$ is an example of an ``important" number in this sense) or as the result of a shared mathematical culture ($10^9$  is more important than $9^{10}$  because we use a decimal notation). The concept of algorithmic complexity \cite{kolmo,chaitin,levin} (also known as Kolmogorov or Kolmogorov-Chaitin complexity) will be used to explain the curve shape as an ``objective" measure. However, the observed curve is not conform to the curve predicted by an analysis based on algorithmic complexity because of a clear gap separating the distribution into two clouds of points. A clear zone in the value of $N(n)$ first noticed by Philippe Guglielmetti\footnote{On his site \url{http://drgoulu.com/2009/04/18/nombres-mineralises/} last consulted 1 June, 2011.}. We shall call this gap ``Sloane's gap". 
\end{abstract}

\section{Introduction}

The Sloane encyclopedia of integer sequences\cite{sloane3}\footnote{The encyclopedia is available at: http://oeis.org/, last consulted  26 may, 2011.} (OEIS) is a remarkable database of sequences of integer numbers, carried out methodically and with determination over forty years\cite{cipra}. As for May 27, 2011, the OEIS contained 189,701 integer sequences. Its compilation has involved hundreds of mathematicians, which confers it an air of homogeneity and apparently some general mathematical objectivity--something we will discuss later on.

When plotting $N(n)$ (the number of occurrences of an integer in the OEIS) two main features are evident:\\
(a) Statistical regression shows that the points $N(n)$ cluster around $k/n^{1.33}$, where $k = 2.53 \times 10^8$.\\
(b) Visual inspection of the graph shows that actually there are two distinct sub-clusters  (the upper one and the lower one) and there is a visible gap between them. We introduce and explain the phenomenon of ``Sloane's gap.''

The paper and rationale of our explanation proceeds as follows:\\
We explain that (a) can be understood using algorithmic information theory. If $U$ is a universal Turing machine, and we denote $m(x)$ the probability that $U$ produces a string $x$, then $m(x) = k 2^{-K(x)+O(1)}$, for some constant $k$, where $K(x)$ is the length of the shortest description of $x$ via $U$. $m(.)$ is usually refered to as the Levin's universal distribution or the Solomonoff-Levin measure \cite{levin}. For a number $n$, viewed as a binary string via its binary representation, $K(n) \leq \log_{2} n + 2 \log_{2} \log_{2} n + O(1)$ and, for most $n$, $K(n) \geq \log_{2} n$. Therefore for most $n$, $m(n)$ lies between $k/( n (\log_{2} n)^2)$ and $k/n$. Thus, if we view OEIS in some sense as a universal Turing machine, algorithmic probability explains (a).

Fact (b), however, is not predicted by algorithmic complexity and is not produced when a database is populated with automatically generated sequences. This gap is unexpected and requires an explanation. We speculate that OEIS is biased towards social preferences of mathematicians and their strong interest in certain sequences of integers (even numbers, primes, and so on). We quantified such a bias and provided statistical facts about it.

\section{Presentation of the database}

The encyclopedia is represented as a catalogue of sequences of whole numbers and not as a list of numbers.  However, the underlying vision of the work as well as its arrangement make it effectively a dictionary of numbers, with the capacity to determine the particular properties of a given integer as well as how many known properties a given integer possesses. 

A common use of the Sloane encyclopedia is in determining the logic of a sequence of integers. If, for example, you submit to it the sequence 3, 4, 6, 8, 12,  14,  18,  20..., you will instantly find that it has to do with the sequence of prime augmented numbers, as follows: 2+1, 3+1, 5+1, 7+1, 11+1, 13+1, 17+1, 19+1...

Even more interesting, perhaps, is the program's capacity to query the database about an isolated number. Let us take as an example the  Hardy-Ramanujan number, 1729 (the smallest integer being the sum of two cubes of two different shapes). The program indicates that it knows of  more than 350 sequences to which 1729 belongs. Each one identifies a property of 1729 that it is possible to examine. The responses are classified in order of importance, an order based on the citations of sequences in mathematical commentaries and the encyclopedia's own cross-references. Its foremost property is that it is the third  Carmichael number (number $n$ not prime for which $\forall a\in \mathbb{N}^{\ast },$ $n|a^{n}-a)$). Next in importance is that 1729 is the sixth pseudo prime in base 2 (number $n$ not prime such that $n|2^{n-1}-1$).  Its third property is that it belongs among the terms of a simple generative series. The property expounded by Ramanujan from his hospital bed appears as the fourth principle. In reviewing the responses from the encyclopedia, one finds further that:
\begin{itemize}
\item 1729 is the thirteenth number of the form $n^3+1$;
\item 1729 is the fourth ``factorial sextuple", that is to say, a product of successive terms of the form $6n +1$: $1729=1\times 7\times 13\times 19$;
\item 1729 is the ninth number of the form $n^{3}+(n+1)^{3}$;
\item 1729 is the sum of the factors of a perfect square ($33^2$);
\item 1729 is a number whose digits, when added together yield its largest factor ($1+7+2+9=19$ and $1729 =7\times 13\times 19$);
\item 1729 is the product of 19 a prime number, multiplied by 91, its inverse;
\item 1729 is the total number of ways to express 33 as the sum of 6 integers.
\end{itemize}

The sequence encyclopedia of Neil Sloane comprises more than 150\,000 sequences.  A partial version retaining only the most important sequences of the database was published by Neil Sloane and Simon Plouffe\cite{sloane2} in 1995. It records a selection of 5487 sequences\cite{sloane2} and echoes an earlier publication by Sloane \cite{sloane}.

\begin{figure}[htdp]
\centering
   \scalebox{.6}{\includegraphics{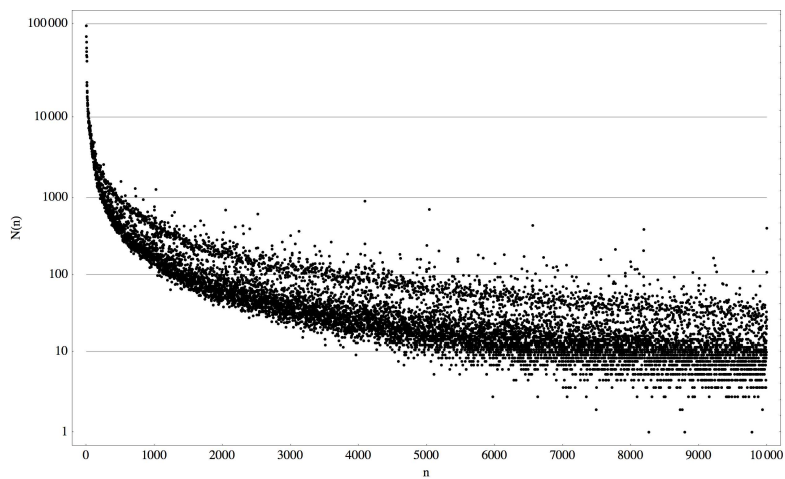}}
\caption{Number of occurrences of $N(n)$ as a function of $n$ per $n$ ranging from 1 to 10\,000. Logarithmic scale in ordinate.}
\end{figure}

Approximately forty mathematicians constitute the ``editorial committee" of the database, but any user may propose sequences. If approved, they are added to the database according to criteria of mathematical interest. Neil Sloane's flexibility is apparent in the  ease with which he adds new sequences as they are proposed. A degree of filtering is inevitable to maintain the quality of the database. Further, there exist a large number of infinite families of sequences (all the sequences of the form ($kn$), all the sequences of the form ($k^n$), etc.), of which it is understood that only the first numbers are recorded in the encyclopedia. A program is also used in the event of a failure of a direct query which allows sequences of families that are not explicitly recorded in the encyclopedia to be recognized. 

Each sequence recorded in the database appears in the form of its first terms. The size of first terms associated with each sequence is limited to approximately 180 digits. As a result, even if the sequence is easy to calculate, only its first terms will be expressed. Next to the first terms and extending from the beginning of the sequence, the encyclopedia proposes all sorts of other data about the sequence, e.g., the definitions of it and bibliographical references. 

Sloane's integer encyclopedia is available in the form of a data file that is easy to read, and that contains only the terms retained for each sequence. One can download the data file free of charge and use it--with mathematical software, for example--to study the expressed numbers and conduct statistical research about the givens it contains.

One can, for example, ask the question: ``Which numbers do not appear in Sloane's encyclopedia?" At the time of an initial calculation conducted in August 2008 by Philippe Guglielmetti, the smallest absent number tracked down was 8795, followed in order by 9935, 11147, 11446, 11612, 11630,... When the same calculation was made again in February 2009, the encyclopedia having been augmented by the addition of several hundreds of new sequences, the series of absent numbers was found to comprise 11630, 12067, 12407, 12887, 13258...

The instability over time of the sequence of missing numbers in the OEIS suggests the need for a study of the distribution of numbers rather than of their mere presence or absence. Let us consider the number of properties of an integer, $N(n)$, while measuring it by the number of times $n$ appears in the number file of the Sloane encyclopedia. The sequence $N(n)$ is certainly unstable over time, but it varies slowly, and certain ideas that one can derive from the values of $N(n)$ are nevertheless quite stable. The values of $N(n)$ are represented in Figure 1. In this logarithmic scale graph a cloud formation with regular decline curve is shown.

Let us give a few examples: the value of $N(1729)$ is 380 (February 2009), which is fairly high for a number of this order of magnitude. For its predecessor, one nevertheless calculates $N(1728)=622$, which is better still.  The number 1728 would thus have been easier for Ramanujan!  Conversely, $N(1730)=106$ and thus 1730 would have required a more elaborate answer than 1729.  

The sequence $\left( N(n)\right) _{n\in \mathbb{N}^{\ast }}$ is generally characterized by a decreasing  curve. However, certain numbers $n$ contradict this rule and possess more properties than their predecessors: $N(n)>N(n-1)$.

We can designate such numbers as ``interesting". The first interesting number according to this definition is 15, because $N(15)=34\,183$ and $N(14)=32\,487$. Appearing next in order are 16, 23, 24, 27, 28, 29, 30, 35, 36, 40, 42, 45, 47, 48, 52, 53, etc. 

We insist on the fact that, although unquestionably dependent on certain individual decisions made by those who participate in building the sequence database, the database is not in itself arbitrary. The number of contributors is very large, and the idea that the database represents an objective view (or at least an intersubjective view) of the numeric world could be defended on the grounds that it comprises the independent view of each person who contributes to it and reflects a stable mathematical (or cultural) reality.

Indirect support for the idea that the encyclopedia is not arbitrary, based as it is on the cumulative work of the mathematical community, is the general cloud-shaped formation of points determined by $N(n)$, which aggregates along a regular curve (see below).

Philippe Guglielmetti has observed that this cloud possesses a remarkable characteristic\footnote{Personal communication with one of the authors, 16th of February, 2009.}: it is divided into two parts separated by a clear zone, as if the numbers sorted themselves into two categories, the more interesting above the clear zone, and the less interesting below the clear zone.  We have given the name ``Sloane's Gap" to the clear zone that divides in two the cloud representing the graph of the function $n\longmapsto N(n)$. 
Our goal in this paper is to describe the form of the cloud, and then to formulate an explanatory hypothesis for it. 

\section{Description of the cloud}
Having briefly described the general form 
of the cloud, we shall direct ourselves more particularly to the gap, and we will investigate what characterizes the points that are situated above it.

\subsection{General shape}

The number of occurrences $N$ is close to a grossly decreasing convex function of $n$, as one can see from Figure 1.

A logarithmic regression provides a more precise idea of the form of the cloud for $n$ varying from 1 to 10\,000. In this interval, the coefficient of determination of the logarithmic regression of $\ln \left( N\left( n\right) \right)$ in $n$ is of $r^2=.81$, and the equation of regression gives the estimation:
$$\ln \left( N\left( n\right) \right) \simeq -1.33\ln (n)+14.76$$
or
$$\hat{N}\left( n\right) =\frac{k}{n^{1.33}},$$ 
where $k$ is a constant having the  approximate value $2.57\times 10^8$, and $\hat{N}$ is the estimated value for $N$.

Thus the form of the function $N$ is determined by the equation above.  Is the existence of Sloane's gap natural then, or does it demand a specific explanation? We note that to our knowledge, only one publication mentions the existence of this split \cite{delahaye}.

\begin{figure}[htdp]
\centering
   \scalebox{.5}{\includegraphics{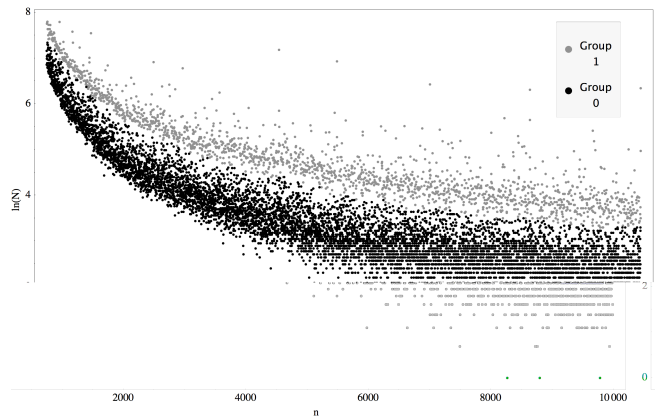}}
\caption{The curve represents the logarithmic regression of $\ln(N)$ as a function of $n$ for $n$ varying from 1 to 10\,000. The grey scale points are those that are classified as being ``above" the gap, while the others are classified as being ``below" it.}
\end{figure}

\subsection{Defining the gap}
In order to study the gap, the first step is to determine a criterion for classification of the points. Given that the ``gap" is not clearly visible for the first values of $n$, we exclude from our study numbers less than 300. 

One empirical method of determining the boundary of the gap is the following: for the values ranging from 301 to 499, we use a straight line adjusted ``by sight", starting from the representation of $\ln(N)$ in functions of $n$. For subsequent values, we take as limit value of $n$ the 82nd percentile of the interval $[n-c,n+c]$. $c$ is fixed at 100 up to $n=1000$, then to 350. It is clearly a matter of a purely empirical choice that does not require the force of a demonstration.  The result corresponds roughly to what we perceive as the gap, with the understanding that a zone of uncertainty will always exist, since the gap is not entirely devoid of points. Figure 2 shows the resulting image.

\subsection{Characteristics of numbers ``above"}
We will henceforth designate as $A$ the set of abscissae of points classified ``above" the gap by the method that we have used. Of the numbers between 301 and 10\,000, 18.2\%  are found in A-- 1767 values.

In this section, we are looking for the properties of these numbers.  Philippe Guglielmetti has already remarked that the prime numbers and the powers of two seem to situate themselves more frequently above the gap.  The idea is that certain classes of numbers that are particularly simple or of particular interest to the mathematician are over-represented.

\subsubsection{Squares}
83 square numbers are found between 301 and 10\,000. Among these, 79 are located above the gap, and 4 below the gap, namely, numbers 361, 484, 529, and 676.  Although they may not be elements of $A$, these numbers are close to the boundary.  One can verify that they collectively realize the local maximums for $\ln (N)$ in the set of numbers classified under the cloud.  One has, for example, $N(361)=1376$, which is the local maximum of $\left\{ N\left( n\right), n\in \left[ 325,10~000\right] \backslash A\right\} $. For each of these four numbers, Table 1 gives the number of occurrences $N$ in Sloane's list, as well as the value limit that they would have to attain to belong to $A$.

95.2\% of squares are found in $A$, as opposed to 17.6\% of non-squares.  The probability that a square number will be in $A$ is thus 5.4 times greater than that for the other numbers.

\begin{equation*}
\begin{tabular}{c|c|c}
\hline
$n$ & $N\left( n\right) $ & value limit \\ 
\hline
361 & 1376 & 1481 \\ 
484 & 976 & 1225 \\ 
529 & 962 & 1065 \\ 
676 & 706 & 855\\
\hline
\end{tabular}
\end{equation*}

\begin{center}
Table 1--List of the square numbers $n$ found between 301 and 10\,000 not belonging to $A$, together with their frequency of occurrence and the value of $N(n)$ needed for $n$ to be classified in $A$.
\end{center}

\subsubsection{Prime numbers}
The interval under consideration contains 1167 prime numbers. Among them, 3 are not in $A$: the numbers 947, 8963, and 9623.  These three numbers are very close to the boundary. 947 appears 583 times, while the limit of $A$ is 584. Numbers 8963 and 6923 appear 27 times each, and the common limit is 28. 

99.7\% of prime numbers belong to $A$, and 92.9\% of non-prime numbers belong to the complement of $A$.  The probability that a prime number will belong to $A$ is thus 14 times greater than the same probability for a non-prime number.

\subsubsection{A multitude of factors}
Another class of numbers that is seemingly over-represented in set $A$ is the set of integers that have ``a multitude of factors". This  is based on the observation that the probability of belonging to $A$ increases with the number of prime factors (counted with their multiples), as can be seen in Figure 3. To refine this idea we have selected the numbers $n$ of which the number of prime factors (with their multiplicty) exceeds the 95th percentile, corresponding to the interval $[n-100,n+100]$.

811 numbers meet this criterion. Of these, 39\% are found in $A$, as opposed to 16.3\% for the other numbers.  The probability that a number that has a multitude of prime factors will belong to $A$ is thus 2.4 times greater than the same probability for a number that has a smaller number of factors. Table 2 shows the composition of $A$ as a function of the classes that we have considered.

\begin{figure}[htdp]
\centering
   \scalebox{.7}{\includegraphics{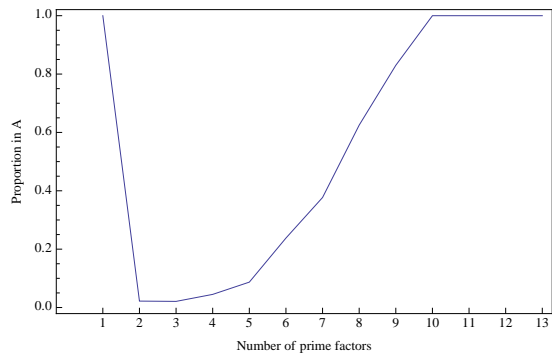}}
\caption{For each number of prime factors (counted with their multiples) one presents, the proportion of integers belonging to $A$ is given.  For the interval determined above, all numbers with at least 10 factors are in $A$.}
\end{figure}

\begin{equation*}
\begin{tabular}{c|c|c|c}
\hline
class & number in $A$ & \% of $A$ & \% (cumulated) \\ 
\hline
primes & 1164 & 65.9 & 65.9 \\ 
squares & 79 & 4.5 & 70.4 \\ 
many factors & 316 & 17.9 & 87.9\\
\hline
\end{tabular}
\end{equation*}

\begin{center}
Table 2--For each class of numbers discussed above, they give the number of occurrences in $A$, the corresponding percentage and the cumulative\\ percentage in $A$.
\end{center}

\subsubsection{Other cases}
The set $A$ thus contains almost all prime numbers, 95\% of squares, and a significant percentage of numbers that have a multitude of factors and all the numbers possessing at least ten prime factors (counted with their multiplicity).

These different classes of numbers by themselves represent 87.9\%  of $A$. Among the remaining numbers, some evince outstanding properties, for example, linked to decimal notation, as in: 1111, 2222, 3333$\ldots$. Others have a simple form, such as 1023, 1025, 2047, 2049... that are written $2^n+1$ or $2^n-1$.

When these cases that for one reason or another possess an evident ``simplicity" are eliminated, there remains a proportion of less than 10\%  of numbers in $A$ for which one cannot immediately discern any particular property.

\section{Explanation of the cloud-shape formation}
\subsection{Overview of the theory of algorithmic complexity}Save in a few exceptional cases,
 for a number to possess a multitude of properties implies  that the said properties are simple, where simple is taken to mean ``what may be expressed in a few words". Conversely, if a number possesses a simple property, then  it will  possess many properties.  For example, if $n$ is a multiple of 3, then $n$ will be a even multiple of 3 or a odd multiple of 3. Being a ``even multiple of 3" or ``odd multiple of 3" is a little more complex than just being a ``multiple of 3", but it is still simple enough, and one may further propose that many sequences  in Sloane's database are  actually sub-sequences of other, simpler ones.  In specifying a simple property, its definition becomes more complex (by generating a sub-sequence of itself), but since there are many ways to specify a simple property, any number that possesses a simple property necessarily possesses numerous properties that are also simple. 

The property of $n$ corresponding to a high value of $N(n)$ thus seems related to the property of admitting a ``simple" description.  The value $N(n)$ appears in this context as an indirect measure of the simplicity of $n$, if one designates as  ``simple" the numbers that have properties expressible in a few words. 

Algorithmic complexity theory\cite{kolmo,chaitin,levin} assigns a specific mathematical sense to the notion of simplicity, as the objects that ``can be described with a short definition". Its modern formulation can be found in the work of Li and Vitanyi\cite{li}, and Calude\cite{calude}.

Briefly, this theory proposes to measure the complexity of a finite object in binary code (for example, a number written in binary notation) by the length of the shortest program that generates a representation of it. The reference to a universal programming language (insofar as all computable functions can possess a program) leads to a theorem of invariance that warrants a certain independence of the programming language. 

More precisely, if $L_{1}$ and $L_{2}$ are two universal languages, and if one notes $K_{L_{1}}$ (resp. $K_{L_{2}}$) algorithmic complexity defined with reference to $L_{1}$ (resp. to $L_{2}$), then there exists a constant $c$ such that $|K_{L_1}(s)-K_{L_2}(s)| < c $ for all finite binary sequences $s$.

A theorem (see for example [theorem 4.3.3. page 253 in \cite{li}]) links the probability of obtaining an object $s$ (by activating a certain type of universal TM--called optimal--running on binary input where the bits are chosen uniformly random) and its complexity $K(s)$. The rationale of this theorem is that if a number has many properties then it also has a simple property.

The translation of this theorem for $N(n)$ is that if one established a universal language $L$, and established a complexity limit $M$ (only admitting descriptions of numbers capable of expression in fewer than $M$ symbols), and counted the number of descriptions of each integer, one would find that  $\frac{N(n)}{M}$ (where $M=\sum_{i\in \mathbb{N}}N(i)$) is approximately proportional to: 
$\frac{1}{2^{K\left( n\right) }}$: $$\frac{N(n)}{M}=\frac{1}{2^{K\left( n\right) +O(\ln (\ln (n)) )}}.$$

Given that $K(n)$ is non computable because of the undecidability of the halting problem and the role of the additive constants involved, a precise calculation of the expected value of $N(n)$ is impossible. By contrast, the strong analogy between the theoretical situation envisaged by algorithmic complexity and the situation one finds when one examines $N(n)$ inferred from Sloane's database, leads one to think that $N(n)$ should be asymptotically dependent on $\frac{1}{2^{K\left( n\right) }}$.  Certain properties of $K(n)$ are obliquely  independent of the reference language chosen to define $K$. The most important of these are:

\begin{itemize}
\item $K(n)<\log _{2}(n)+2\log _{2}(\log_{2}(n))+c^{\prime }$ ($c^{\prime }$ a constant)
\item the proportion of $n$ of a given length (when written in binary) for which $K(n)$ recedes from $\log_2(n)$ decreases exponentially (precisely speaking, less than an integer among $2^q$ of length $k$, has an algorithmic complexity  $K(n)\leq k-q$).
\end{itemize}

In graphic terms, these properties indicate that the cloud of points obtained from writing the following $\frac{1}{2^{K\left( n\right) }}$ would be situated above a curve defined by $$f(n)\approx \frac{h}{2^{\log _{2}\left( n\right) }}=\frac{h}{n}$$ ($h$ being a constant), and that all the points  cluster on the curve, with the density of the points deviating from the curve decreasing rapidly.  

This is indeed the situation we observe in examining the curve giving $N(n)$. The theory of algorithmic information thus provides a good description of what is observable from the curve  $N(n)$.  That justifies an a posteriori recourse to the theoretical concepts of algorithmic complexity in order to understand the form of the curve $N(n)$.  By contrast, nothing in the theory leads one to expect a gap like the one actually observed. To the contrary, continuity of form is expected from the fact that $n+1$ is never much more complex than $n$.

To summarize, if $N(n)$ represented an objective measure of the complexity of numbers (the larger $N(n)$ is, the simpler $n$ ), these values would then be comparable to those that yield $\frac{1}{2^{K\left( n\right) }}$. One should thus observe a rapid decrease in size, and a clustering of values near the base against an oblique curve, but one should not observe a gap, which presents itself here as an anomaly.  

To confirm the conclusion that the presence of the gap results from special factors and render it more convincing, we have conducted a numerical experiment.

\begin{figure}[htdp]
\centering
   \scalebox{.6}{\includegraphics{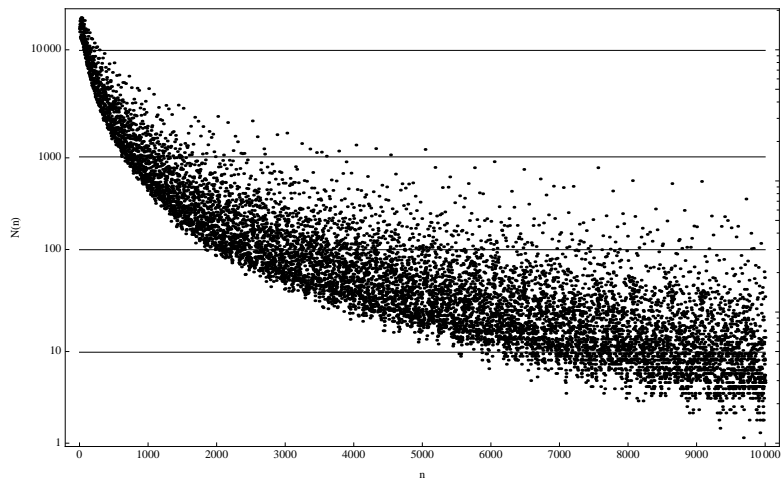}}
\caption{Graph of $N(n)$ obtained with random functions, similar to that belonging to Sloane's Database (Figure 1). Eight million values have been generated.}
\end{figure}

We define random functions $f$ in the following manner (thanks to the algebraic system $Mathematica$): 

\begin{enumerate}
\item Choose at random a number $i$ between 1 and 5 (bearing in mind in the selection the proportions of functions for which $i=1$, $i=2$, $\ldots$, $i =5$ among all those definable in this way).

\item If $i=1$, $f$ is defined by choosing uniformly at random a constant $k\in \left\{ 1,...,9\right\}$, a binary operator $\varphi $ from among the following list: $+$, $\times$, and subtraction sign, in a uniform manner, and a unary operand $g$ that is identity with probability .8, and the function squared with probability .2 (to reproduce the proportions observed in Sloane's database). One therefore posits  $f_i(n)=\varphi (g(n),k).$

\item If $i\geq 2$, $f_i$ is defined by $f_i(n) =\varphi (g(f_{i-1}(n)),k), $ where $k$ is a random integer found between 1 and 9, $g$ and $\varphi $ are selected as described in the point 2 (above), and $f_{i-1}$ is a random function selected in the same manner as in 2.
\end{enumerate}

For each function $f$ that is generated in this way, one calculates $f(n)$ for $n=1, $\ldots$, 20.$  These terms are regrouped and counted as for $N(n)$. The results appear in Figure 4. The result confirms what the relationship with algorithmic complexity would lead us to expect. There is a decreasing oblique curve with a mean near 0, with clustering  of the points near the base, but no gap.

\subsection{The gap: A social effect?}

This anomaly with respect to the theoretical implications and modeling is undoubtedly a sign that what one sees in Sloane's database is not a simple objective measure of complexity (or of intrinsic mathematical interest), but rather a trait of psychological or social origin that mars its pure expression.  That is the hypothesis that we propose here. Under all circumstances, a purely mathematical vision based on algorithmic complexity would encounter an obstacle here, and the social hypothesis is both simple and natural owing to the fact that Sloane's database, while it is entirely ``objective", is also a social construct.

\begin{figure}[htdp]
\centering
   \scalebox{.5}{\includegraphics{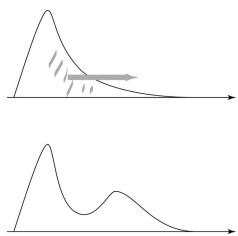}}
\caption{The  top figure above represents the local distribution of $N$ expected without taking into account the social factor.}
\end{figure}

Figure 5 illustrates and specifies our hypothesis that the mathematical community is particularly interested in certain numbers of moderate or weak complexity (in the central zone or on the right side of the distribution), and this interest creates a shift toward the right-hand side of one part of the distribution (schematized here by the grey arrow). The new distribution that develops out of it (represented in the bottom figure) presents a gap.

We suppose that the distribution anticipated by considerations of complexity is deformed by the social effect concomitant with it: mathematicians are more interested in certain numbers that are linked to selected properties by the scientific community.  This interest can have cultural reasons or mathematical reasons (as per results already obtained), but in either case it brings with it  an over-investment on the part of the mathematical community.  The numbers that receive this specific over-investment are not in general complex, since  interest is directed toward them  because certain regularities have been discovered in them. Rather, these numbers are situated near the pinnacle of a theoretical asymmetrical distribution.  Owing to the community's over-investment, they are found to have shifted towards the right-hand side of the distribution, thus explaining Sloane's gap.

It is, for example, what is generated by  numbers of the form $2^n +1$, all in A, because  arithmetical results can be obtained from this type of number that are useful to prime numbers.  Following some interesting preliminary discoveries, scientific investment in this class of integers has become intense, and they appear in numerous sequences.  Certainly, $2^n+1$ is objectively a simple number, and thus it is normal that it falls above the gap.  Nevertheless, the difference in complexity between $2^n+1$ and $2^n+2$ is weak. We suppose that the observed difference also reflects a social dynamic which tends to augment $N(2^n+1)$ for reasons that complexity alone would not entirely explain.

\section{Conclusion}

The cloud of points representing the function $N$ presents a general form evoking an underlying function characterized by rapid decrease and ``clustering near the base" (local asymmetrical distribution). This form is explained, at least qualitatively, by the theory of algorithmic information.

If the general cloud formation was anticipated, the presence of Sloane's gap has, by contrast, proved more challenging to its observers.  This gap has not, to our knowledge, been successfully explained on the basis of uniquely numerical considerations that are independent of human nature as it impinges on the work of mathematics. Algorithmic complexity anticipates a certain ``continuity"  of $N$, since the complexity of $n+1$ is always close to that of $n$. The discontinuity that is manifest in Sloane's gap is thus difficult to attribute to purely mathematical properties independent of social contingencies. 

By contrast, as we have seen, it is explained very well by the conduct of research that entails the over-representation of certain numbers of weak or medium complexity.  Thus the cloud of points representing the function $N$ shows features that can be understood as being the result of at the same time human and purely mathematical factors.

\end{document}